\documentclass[a4paper,12pt,twoside,leqno,final]{amsart}
\usepackage{amsmath}
\usepackage{amssymb}

\newtheorem{thm}{Theorem}[section]
\newtheorem{lem}[thm]{Lemma}

\newtheorem{prop}[thm]{Proposition}
\newtheorem{cor}[thm]{Corollary}

\newcommand{\C}{{\mathbb C}}
\newcommand{\D}{{\mathbb D}}
\newcommand{\R}{{\mathbb R}}
\newcommand{\T}{{\mathbb T}}
\newcommand{\Z}{{\mathbb Z}}
\newcommand{\N}{{\mathbb N}}

\newcommand{\bmo}{{\rm BMO}}
\newcommand{\bmoa}{{\rm BMOA}}

\newcommand{\Al}{A^\alpha}

\newcommand{\La}{\Lambda}

\newcommand{\eps}{\varepsilon}

\newcommand{\f}{\frac}
\newcommand{\ov}{\overline}
\newcommand{\al}{\alpha}
\newcommand{\be}{\beta}
\newcommand{\Ga}{\Gamma}
\newcommand{\ga}{\gamma}

\newcommand{\de}{\delta}

\newcommand{\ze}{\zeta}
\renewcommand{\th}{\theta}
\newcommand{\si}{\sigma}
\newcommand{\ph}{\varphi}

\newcommand{\Om}{\Omega}
\newcommand{\Omte}{\Om(\th,\eps)}

\newcommand{\const}{\text{\rm const}}

\numberwithin{equation}{section}

\title[Factorization theorems for pseudocontinuable functions]
{Factorization and non-factorization theorems\\ 
for pseudocontinuable functions}
\author{Konstantin M. Dyakonov}
\address{Departament de Matem\`atiques i Inform\`atica, IMUB, BGSMath, Universitat de Barcelona, Gran Via 585, E-08007 Barcelona, Spain}
\address{ICREA, Pg. Llu\'is Companys 23, E-08010 Barcelona, Spain}
\email{konstantin.dyakonov@icrea.cat}
\thanks{Supported in part by grant MTM2014-51834-P from El Ministerio de Econom\'ia 
y Competitividad (Spain) and grant 2014-SGR-289 from AGAUR (Generalitat de Catalunya).}

\begin{document}
\begin{abstract}
Let $\theta$ be an inner function on the unit disk, and let $K^p_\theta:=H^p\cap\th\ov{H^p_0}$ be the associated star-invariant subspace of the Hardy space $H^p$, with $p\ge1$. While a nontrivial function $f\in K^p_\theta$ is never divisible by $\theta$, it may have a factor $h$ which is \lq\lq not too different" from $\theta$ in the sense that the ratio $h/\theta$ (or just the anti-analytic part thereof) is smooth on the circle. In this case, $f$ is shown to have additional integrability and/or smoothness properties, much in the spirit of the Hardy--Littlewood--Sobolev embedding theorem. The appropriate norm estimates are established, and their sharpness is discussed.
\end{abstract}

\maketitle

\section{Introduction and results}

The pseudocontinuable functions in the paper's title are the noncyclic vectors of the backward shift operator 
$$S^*:\,f\mapsto\f{f-f(0)}z$$
acting on the Hardy space $H^2$ -- or, more generally, $H^p$ (see the definition below) -- of the unit disk $\D:=\{z\in\C:|z|<1\}$. As usual, \lq\lq noncyclic" means \lq\lq lying in some proper (closed) invariant subspace", and a well-known result from \cite{DSS} tells us that a function $f$ is noncyclic for $S^*$ if and only if it admits a {\it pseudocontinuation} to $\D^-:=\{z:|z|>1\}$. The latter means that there exists a meromorphic function of bounded characteristic in $\D^-$ whose boundary values agree with $f$ almost everywhere on the circle $\T:=\partial\D$. The other key notion in this paper is {\it smoothness}, a phenomenon that can also be described in terms of a \lq\lq pseudocontinuation", this time understood as {\it pseudoanalytic extension} in the sense of Dyn'kin (see \cite{Dyn}). The general principle, as explained in \cite{Dyn}, is that a function $f$ holomorphic on $\D$ (and, say, continuous up to the boundary) will be smooth on $\T$, in some sense or other, if and only if it extends to $\D^-$ as a $C^1$ function whose Cauchy--Riemann $\ov\partial$-derivative becomes appropriately small -- or does not grow too fast -- near $\T$. 

\par Thus, loosely speaking, we are concerned with the interplay of the two kinds of pseudocontinuation. In this connection, we also mention the survey paper \cite{DCora} which summarizes some of the previous results pertaining to the two topics and discusses the interrelationship between them. 

\par Now let us try and describe the setup more accurately. Recall, first of all, that the Hardy space $H^p=H^p(\D)$ with $0<p<\infty$ is formed by those holomorphic functions $f$ on $\D$ which satisfy 
$$\sup_{0<r<1}\int_\T|f(r\ze)|^pdm(\ze)<\infty$$
(we write $m$ for the normalized arclength measure on $\T$), while $H^\infty=H^\infty(\D)$ stands for the space of bounded holomorphic functions. As is customary, we identify $H^p$ functions with their boundary values (defined almost everywhere on the circle) and treat $H^p$ as a subspace of $L^p=L^p(\T,m)$, endowed with the standard $L^p$-norm $\|\cdot\|_p$. Further, a function $\th\in H^\infty$ is said to be {\it inner} if $|\th|=1$ almost everywhere on $\T$. 

\par Beurling's famous theorem characterizes the invariant subspaces of the (forward) shift operator $S:f\mapsto zf$ in $H^2$ as those of the form $\th H^2$, where $\th$ is an inner function; see, e.g., \cite[Chapter II]{G}. Accordingly, the $S^*$-invariant (or {\it star-invariant}) subspaces of $H^2$ can be written as 
$$H^2\ominus\th H^2=:K^2_\th,$$
with $\th$ inner. These are also known as {\it model subspaces}, due to their role in the Sz.-Nagy--Foia\c{s} operator model; see \cite{N}. 

\par Obviously enough, the only function $f\in K^2_\th$ that is divisible by $\th$ (in $H^2$) is $f\equiv0$. The question we address here is whether, and to what extent, this trivial uniqueness theorem is stable under small -- or not so small -- perturbations of the hypothesis. Assume, for instance, that a function $f\in K^2_\th$ is divisible by a certain (possibly non-inner) $h\in H^\infty$ which is \lq\lq not too different" from $\th$, meaning that the ratio $h/\th$ -- or maybe just its projection onto $L^2\ominus H^2$ -- is suitably smooth on $\T$. Does this provide us with some sort of control over the size of $f$? 

\par More generally, given an inner function $\th$, we introduce the star-invariant subspace $K^p_\th$ in $H^p$ by putting 
\begin{equation}\label{eqn:defkpth}
K^p_\th:=H^p\cap\th\,\ov{H^p_0},\qquad1\le p\le\infty, 
\end{equation}
where $H^p_0:=zH^p$ and the bar denotes complex conjugation. (When $p=2$, the new definition coincides with the old one.) Letting 
$$\widetilde f:=\ov z\ov f\th,$$
we may also rewrite \eqref{eqn:defkpth} in the form 
$$K^p_\th=\{f\in H^p:\,\widetilde f\in H^p\},\qquad1\le p\le\infty.$$
In fact, the transformation $f\mapsto\widetilde f$ is an important antilinear isometry (or involution) that maps $K^p_\th$ onto itself. We shall not deal with $K^p_\th$ spaces in the range $0<p<1$, and in most cases we restrict our attention to $K^1_\th$, the largest among the subspaces of this kind that we do consider. 

\par As before, we observe that no nontrivial function in $K^1_\th$ is divisible by $\th$, and we look for quantitative refinements -- or perturbations -- of this fact. Specifically, assume that a function $f\in K^1_\th$ can be factored as $f=gh$, where $g\in H^p$ (not necessarily with $p\ge1$) and $h$ is an $H^\infty$ function which is \lq\lq not too different" from $\th$, in the same sense as above. Then, as it turns out, $f$ is likely to enjoy improved integrability properties, and is thus forced to be smaller than one would {\it a priori} expect. In particular, we shall see that if the ratio $h/\th$ (or just its anti-analytic part) satisfies the Lipschitz or Zygmund condition of order $\al$, with $\al<p^{-1}$, then $f$ will be in $H^r$ with $r=p/(1-\al p)$, this last exponent being sharp. 

\par But first we need to introduce a bit of notation. We let $P_+$ and $P_-$ denote the {\it Riesz projections}, whose action on a function $F\in L^1$ is given by 
$$(P_+F)(z):=\sum_{k\ge0}\widehat F(k)z^k\quad\text{\rm and}\quad (P_-F)(z):=\sum_{k<0}\widehat F(k)z^k;$$ 
here $\widehat F(k):=\int_\T F(\ze)\ov\ze^kdm(\ze)$ is the $k$th Fourier coefficient of $F$. When restricted to $L^p$ with $1<p<\infty$, these become bounded projections onto $H^p$ and onto $\ov{H^p_0}$, respectively; see \cite[Chapter III]{G}. 

\par Further, given a number $\al\in\R$, we denote by $\Al$ the set of those holomorphic functions $f$ on $\D$ which satisfy 
\begin{equation}\label{eqn:nalphaf}
\mathcal N_\al(f):=\sup_{z\in\D}\left|f^{(n)}(z)\right|(1-|z|)^{n-\al}<\infty,
\end{equation}
where $n$ is the least nonnegative integer in the interval $(\al,\infty)$. We then put 
$$\|f\|_{\Al}:=\sum_{k=0}^{n-1}\left|f^{(k)}(0)\right|+\mathcal N_\al(f).$$
When $\al=-\be<0$, one takes $n=0$, so that 
$$\|f\|_{A^{-\be}}=\mathcal N_{-\be}(f)=\sup_{z\in\D}|f(z)|(1-|z|)^\be;$$
the corresponding growth spaces $A^{-\be}$, with $\be>0$, are known as {\it Korenblum classes}. We note that, for $0<p<\infty$, $A^{-1/p}$ contains the Hardy space $H^p$, and moreover, $\|f\|_{A^{-1/p}}\le\|f\|_p$ for all $f\in H^p$. 
\par When $\al=0$, we get the so-called {\it Bloch space}, denoted usually by $\mathcal B$ rather than $A^0$; it will not reappear in what follows. 
\par Finally, when $\al>0$, \eqref{eqn:nalphaf} becomes a smoothness condition on $\T$, and the $\Al$ spaces that arise are the classical {\it analytic Lipschitz--Zygmund spaces}. Equivalently, we have $\Al=\La^\al\cap H^\infty$, where $\La^\al=\La^\al(\T)$ is the appropriate Lipschitz--Zygmund class on the circle. The latter is defined, for $0<\al<\infty$, as the set of those functions $f\in C(\T)$ which obey the condition 
$$\left\|\Delta^n_hf\right\|_\infty=O(|h|^\al),\qquad h\in\R,$$
for some (and then each) integer $n$ with $n>\al$;  
here $\|\cdot\|_\infty$ is the $\sup$-norm on $\T$, and $\Delta_h^n$ denotes the $n$th order difference operator with step $h$. (As usual, the difference 
operators $\Delta_h^k$ are defined by induction: one sets $(\Delta_h^1f)(\ze):=f(e^{ih}\ze)-f(\ze)$ and 
$\Delta_h^kf:=\Delta_h^1\Delta_h^{k-1}f$.) Furthermore, decomposing a given function $f\in\La^\al$ as $f=\varphi+\ov\psi$ with $\varphi\in H^2$ and $\psi\in H^2_0$, we know from the Privalov--Zygmund theorem that both $\varphi(=P_+f)$ and $\psi(=\ov{P_-f})$ will be in $\Al$, and we norm $\La^\al$ by putting 
$$\|f\|_{\La^\al}:=\|\varphi\|_{\Al}+\|\psi\|_{\Al}.$$

\par Now we are in a position to state our first result. 

\begin{thm}\label{thm:factalphapr} Let $\al$ and $p$ be positive numbers such that $\al<1/p<\al+1$, and 
let $r=p/(1-\al p)$. Suppose that $\th$ is an inner function and $f\in K^1_\th$. Assume, finally, that $f=gh$, 
where $g\in H^p$, $h\in H^\infty$ and $P_-(h\ov\th)\in\La^\al$. Then $f\in H^r$ and 
\begin{equation}\label{eqn:normalphapr}
\|f\|_r\le C\|g\|_p\left\|P_-(h\ov\th)\right\|_{\La^\al},
\end{equation}
with some constant $C=C(\al,p)>0$. 
\end{thm}

One of the two inequalities imposed on the parameters, namely $1/p<\al+1$, is not really restrictive. Indeed, if $1/p\ge\al+1$, then $r\le1$ and the conclusion that $f\in H^r$ is automatic (since $f\in K^1_\th$). It is the restriction $\al<1/p$ that matters, and we shall later see what happens beyond this range. 

\par The following corollary deals with the special case where $h$ is (a divisor of) the inner factor of $f$. On the one hand, the presence of a certain inner factor does not improve integrability properties of a function on $\T$. On the other hand, an improvement is ensured in the current situation by Theorem \ref{thm:factalphapr}. By juxtaposing the two facts, one arrives at the result below (of which a formal proof will also be given). 

\begin{cor}\label{cor:innsomealpha} Given an inner function $\th$, suppose that  $f=FI\in K^1_\th$, where $F\in H^1$ and $I$ is inner. If $P_-(I\ov\th)\in\bigcup_{\al>0}\La^\al$, then $f\in\bigcap_{0<p<\infty}H^p$. 
\end{cor}

More can be said when $\th$ is an {\it interpolating Blaschke product}, an assumption to be made in the next corollary. In other words, the role of $\th$ will now (temporarily) be played by
$$B(z):=\prod_j\f{|z_j|}{z_j}\f{z_j-z}{1-\ov z_jz},\qquad\{z_j\}\subset\D,$$
where the product converges and satisfies $\inf_j|B'(z_j)|(1-|z_j|)>0$. (Equivalently, the zeros $\{z_j\}$ of $B$ form an interpolating sequence for $H^\infty$; see \cite[Chapter VII]{G}.) 
\par In view of classical interpolation theorems, it is no wonder that suitable growth restrictions on the values $f(z_j)$ of a function $f\in K^1_B$ may control the global growth of $f$ near $\T$. For instance, the function's membership in $H^p$, $1<p<\infty$, is easily rephrased in terms of the sequence $\{f(z_j)\}$. It is also known that {\it smallness} conditions on $\{f(z_j)\}$ account for {\it smoothness} properties of $\widetilde f:=\ov z\ov fB$, and vice versa; see \cite{DSpb93, DIUMJ94}. The new -- perhaps somewhat surprising -- feature in our next result is that a smallness condition on $\{I(z_j)\}$, where $I$ is the inner factor of $f$, may have a similar effect. 

\begin{cor}\label{cor:smallinnfact} Let $0<\al<\infty$ and let $B$ be an interpolating Blaschke product with zeros $\{z_j\}_{j=1}^\infty$. Assume also that $f=FI\in K^1_B$, where $F$ is outer, $I$ is inner, and  
\begin{equation}\label{eqn:smallatzeros}
\sup_j\f{|I(z_j)|}{(1-|z_j|)^\al}<\infty.
\end{equation}
Then $\widetilde f:=\ov z\ov fB$ is in $\Al$, and so is $F$. 
\end{cor}

Thus, in particular, if the values $I(z_j)$ decay at a certain rate, then $F$ is forced to be smooth (and hence, {\it a fortiori}, bounded). An extreme case of this phenomenon is illustrated by the following example. Let $B_1$ be the Blaschke product with zeros $\{z_j\}_{j=2}^\infty$, so that $B_1=B/b_1$, where $b_1$ is the elementary Blaschke factor with zero at $z_1$. Now, for a function $F\in H^1$, we have $f:=FB_1\in K^1_B$ if and only if $F\in K^1_{b_1}$, and this means that $F$ (as well as $\widetilde f$) is a constant multiple of $(1-\ov z_1z)^{-1}$. The hypotheses and the conclusion of Corollary \ref{cor:smallinnfact} are in this case obviously fulfilled, with $I=B_1$, for each $\al>0$. 

\par Less trivial examples are constructed by perturbing the above picture. Namely, we may replace $B_1$ by a Blaschke product $\mathcal B_1$ with zeros $\ze_j$ ($j=2,3,\dots$), which are taken to be suitably close to the respective $z_j$'s. In fact, it is known (see \cite[Theorem 6.10]{AKS}) that $\mathcal B_1$ will be the inner factor of some non-null $f$ in $K_B^2$ whenever 
\begin{equation}\label{eqn:zerosaks}
\sum_{j=2}^\infty\f{|\ze_j-z_j|}{1-|z_j|}<\infty.
\end{equation}
The function $F:=f/\mathcal B_1$ will then be outer. Also, one easily verifies that 
$$|\mathcal B_1(z_j)|\le\f{|\ze_j-z_j|}{1-|z_j|}\qquad(j=2,3,\dots).$$ 
Therefore, given $\al>0$, we can make the distances $|\ze_j-z_j|$ appropriately small to ensure that both \eqref{eqn:zerosaks} and \eqref{eqn:smallatzeros} hold true, the latter with $I=\mathcal B_1$. Corollary \ref{cor:smallinnfact} will then again be applicable. 

\par Going back to Theorem \ref{thm:factalphapr}, we supplement it with an alternative version, stated as Theorem \ref{thm:factalphaprbis} below, which covers the endpoint case $\al=1/p$. In addition, the factor $\|g\|_p$ on the right-hand side of \eqref{eqn:normalphapr} will be replaced by a smaller quantity involving 
$$\|g\|_{A^{-1/p}}:=\sup\left\{|g(z)|\cdot(1-|z|)^{1/p}:\,z\in\D\right\},$$
the norm of $g$ in $A^{-1/p}$. (Recall that $H^p\subset A^{-1/p}$, the inclusion map being contractive.) On the other hand, we need to replace the hypothesis $P_-(h\ov\th)\in\La^\al$ of Theorem \ref{thm:factalphapr} by a stronger assumption, which we now describe. 

\par Given an inner function $\th$ and a number $\al>0$, a function $h\in H^\infty$ will be called {\it $(\al,\th)$-nice} if $h\th^k\in\La^\al$ for all $k\in\Z$. (Of course, $h$ itself must then be in $\Al$.) A characterization of such functions $h$ can be found in \cite{DSpb93} (see also \cite{DAJM} and \cite{DActa} for related results), and we shall reproduce it as Lemma \ref{lem:kdmultdiv} below. The lemma will tell us, in particular, that a necessary and sufficient condition for $h\in\Al$ to be $(\al,\th)$-nice is that $h\ov\th^N\in\La^\al$ for some (any) fixed integer $N$ exceeding $\al$. The appropriate quantitative characteristic of $(\al,\th)$-niceness may therefore be introduced, for $h\in H^\infty$, by putting
$$\|h\|_{\al,\th}:=\|h\|_{\La^\al}+\|h\ov\th^n\|_{\La^\al}$$
with $n=[\al]+1$ (as usual, we write $[\al]$ for the integral part of $\al$). 

\begin{thm}\label{thm:factalphaprbis} Let $\al$ and $p$ be positive numbers satisfying $\al\le 1/p<\al+1$, and 
let $r=p/(1-\al p)$. Suppose that $\th$ is an inner function and $f\in K^1_\th$. Assume also that $f=gh$, where $g\in H^p$ and $h$ is $(\al,\th)$-nice. 
Then $f\in H^r$ and 
\begin{equation}\label{eqn:normalphaprbis}
\|f\|_r\le C\|g\|^{\al p}_{A^{-1/p}}\|g\|^{1-\al p}_p\|h\|_{\al,\th},
\end{equation}
where $C=C(\al,p)$ is a positive constant. In particular, when $\al p=1$, we have 
$f\in H^\infty$ and 
\begin{equation}\label{eqn:normalphainfty}
\|f\|_\infty\le\const\cdot\|g\|_{A^{-\al}}\|h\|_{\al,\th},
\end{equation}
where the constant depends only on $\al$.
\end{thm}

\par A natural question that arises in connection with Theorems \ref{thm:factalphapr} and \ref{thm:factalphaprbis} is whether the \lq\lq critical exponent" $r=p/(1-\al p)$ is actually sharp. And if it is, one might still ask whether it could be improved under additional smoothness assumptions on $\th$, 
e.g., in the case where $\th'\in H^\ga$ for some $\ga\in(0,1)$. After all, such hypotheses on $\th$ are known to imply certain Sobolev- or Besov-type properties for functions in the associated star-invariant subspace (see \cite{C0, DIUMJ94}), to say nothing of the extreme case $\ga=1$ when the functions involved become analytic on $\T$. Our next proposition shows that the current value of $r$ in the above theorems is indeed optimal, even for very good $\th$'s. 

\begin{prop}\label{prop:sharpness} There is an inner function $\th$ with $\th'\in\bigcap_{0<\ga<1}H^\ga$ that has the following property: whenever 
$\al$, $p$ and $\si$ are positive numbers satisfying $\al<1/p<\al+1$ and $\si>p/(1-\al p)$, one can find 
a function $f\in K^1_\th\setminus H^\si$ representable as $f=gh$, where $g\in H^p$ and $h$ is $(\al,\th)$-nice. 
\end{prop}

\par In fact, the proof will reveal that if $\th$ is chosen appropriately, and if $s>r:=p/(1-\al p)$, then {\it every} $f\in H^s$ (and hence, in particular, every $f\in K^s_\th\setminus H^\si$ with $\si>s$) is of the form $f=gh$, where $g\in H^p$ and $h$ is $(\al,\th)$-nice. On the other hand, Theorem \ref{thm:factalphaprbis} tells us that {\it no} function $f\in K^1_\th\setminus H^r$ admits such a factorization. Thus, we are faced with an amusing \lq\lq all or nothing" dichotomy. 

\par One association that comes to mind in connection with Theorem \ref{thm:factalphapr} is the Hardy--Littlewood--Sobolev embedding theorem (see \cite{HL}): if $0<\al<1/p$ and if $f$ is an analytic function on $\D$ whose (possibly fractional) derivative $f^{(\al)}$ of order $\al$ is in $H^p$, then $f\in H^r$ with $r=p/(1-\al p)$. Furthermore, for $\al>1/p$, the assumption that $f^{(\al)}\in H^p$ forces $f$ to be in $A^\be$ with $\be=\al-1/p$. 

\par In view of this vague analogy, one might try to guess what happens to Theorem \ref{thm:factalphapr} in the latter case (i.e., when $\al>1/p$). A naive conjecture would probably be that in this range, the conclusion that $f\in H^r$ should again be replaced by $f\in A^\be$ with $\be=\al-1/p$. However, no such thing is actually true. To see why, suppose that $\th$ is a nonrational inner function with $\th(0)=0$. The function $f=\th/z$, which belongs to $K^\infty_\th(\subset K^1_\th)$, can be written as $gh$, where $g=1$ and $h=f$. Obviously, $g$ is in $H^\infty$, and hence in every $H^p$, while $P_-(h\overline\theta)=\overline z$ is in $\La^\al$ for every $\al>0$. At the same time, $f$ does not belong to any $A^\be$ space with $\be>0$; it is not even continuous on $\T$. 

\par What {\it is} true, however, in the range $\al>1/p$ is that the \lq\lq partner" 
$\widetilde f$ of $f\in K^1_\th$ must be in $A^\be$, with $\be=\al-1/p$, provided that $f$ can be factored as in Theorem \ref{thm:factalphapr}. In the borderline case $\al=1/p$, this remains true if $A^\be$ gets replaced by $\bmoa$, the analytic subspace of $\bmo=\bmo(\T)$; see \cite[Chapter VI]{G} for the definition and basic properties of the latter. Recalling that $\bmoa=(H^1)^*$, under the usual pairing, we equip $\bmoa$ with the dual space norm $\|\cdot\|_*$. 

\begin{thm}\label{thm:embedlip} Let $\al$ and $p$ be positive numbers with $\al\ge1/p$. Suppose that $\th$ is an inner function and $f\in K^1_\th$. Assume, finally, that $f=gh$, where $g\in H^p$, $h\in H^\infty$ and $P_-(h\ov\th)\in\La^\al$. 
\par{\rm (A)} If $\al>1/p$, then $\widetilde f:=\ov z\ov f\th$ is in $A^\be$, where $\be=\al-1/p$, and 
\begin{equation}\label{eqn:normalphabeta}
\|\widetilde f\|_{\La^\be}\le C\|g\|_p\left\|P_-(h\ov\th)\right\|_{\La^\al}
\end{equation}
with some constant $C=C(\al,p)>0$. 
\par{\rm (B)} If $\al=1/p$, then $\widetilde f\in\bmoa$ and 
\begin{equation}\label{eqn:normalphabmo}
\|\widetilde f\|_*\le C\|g\|_p\left\|P_-(h\ov\th)\right\|_{\La^\al}
\end{equation}
with some $C=C(\al)>0$.
\end{thm}

\par The exponent $\be=\al-1/p$ in part (A) is again sharp, even for nice $\th$'s; this is illustrated by Proposition \ref{prop:sharpnesslip} below. Part (B) should be compared with the $\al=1/p$ case of Theorem \ref{thm:factalphaprbis}, where a stronger assumption allowed us to arrive at the stronger conclusion that $f$ (and hence also $\widetilde f$) is in $H^\infty$. 

\begin{prop}\label{prop:sharpnesslip} There is an inner function $\th$ with $\th'\in\bigcap_{0<\ga<1}H^\ga$ that has the following property: whenever 
$\al$, $p$ and $\de$ are positive numbers satisfying $0<\al-p^{-1}<\de$, one can find a function $f\in K^1_\th$ representable as $f=gh$, with $g\in H^p$, $h\in H^\infty$ and $P_-(h\ov\th)\in\La^\al$, for which $\widetilde f\notin A^\de$.
\end{prop}

\par Our last theorem can be viewed as a variant of Theorem \ref{thm:embedlip}, part (A). The factor $\|g\|_p$ on the right-hand side of \eqref{eqn:normalphabeta} will now be reduced to $\|g\|_{A^{-1/p}}$, while the other factor, $\left\|P_-(h\ov\th)\right\|_{\La^\al}$, will be replaced by a larger quantity. 

\begin{thm}\label{thm:embedlipbis} Let $\al$ and $p$ be positive numbers with $\al>1/p$. Suppose that $\th$ is an inner function and $f\in K^1_\th$. Assume also that $f=gh$, where $g\in H^p$ and $h$ is $(\al,\th)$-nice. Then 
$\widetilde f:=\ov z\ov f\th$ is in $A^\be$, where $\be=\al-1/p$, and 
\begin{equation}\label{eqn:normalphabetabis}
\|\widetilde f\|_{\La^\be}\le C\|g\|_{A^{-1/p}}\|h\|_{\al,\th},
\end{equation}
where $C=C(\al,p)$ is a positive constant.
\end{thm}

\par The next section contains some preliminary facts we need to lean upon, and the rest of the paper is devoted to proving our results. Throughout, the following standard notation will be used: given two nonnegative quantities $U$ and $V$, we write $U\lesssim V$ (or equivalently, $V\gtrsim U$) to mean that $U\le\const\cdot V$. When $U\lesssim V\lesssim U$, we write $U\asymp V$ and say that the two quantities are comparable to one another. Typically, the constants involved in such comparison relations will only depend on parameters such as $\al$, $p$, etc. 

\section{Preliminaries}

Several auxiliary results will be needed. The first of these, stated as Lemma \ref{lem:drsduality} below, is the classical Duren--Romberg--Shields theorem (see \cite{DRS}) which allows us to view the analytic Lipschitz space $\Al$ as $(H^p)^*$, the dual of a certain Hardy space $H^p$ with $0<p<1$. 

\begin{lem}\label{lem:drsduality} Let $0<\al<\infty$ and $p=(1+\al)^{-1}$. Every function $f\in\Al$ induces a continuous linear functional on $H^p$ by the rule 
\begin{equation}\label{eqn:drsdual}
g\mapsto\lim_{\rho\to1^-}\int_\T\ov{f(\ze)}g(\rho\ze)\,dm(\ze),\qquad g\in H^p;
\end{equation}
in particular, the limit on the right exists. Conversely, every functional in $(H^p)^*$ arises in this way. Moreover, the norm of $f$ in $\Al$ is comparable to that of the associated functional \eqref{eqn:drsdual} in $(H^p)^*$, with comparison constants depending only on $\al$.
\end{lem}

\par As another preliminary result, we list some facts about the so-called {\it Carleson curves} associated with an inner function; see \cite[Chapter VIII]{G} for a proof. 

\begin{lem}\label{lem:carlesoncurves} Given an inner function $\th$ and a number $\eps\in(0,1)$, there exists a countable 
(possibly finite) system $\Ga_\eps=\Ga_\eps(\th)$ of simple closed rectifiable curves in $\D\cup\T$ with the following properties. 
\par{\rm (a)} The interiors of the curves in $\Ga_\eps$ are pairwise disjoint; the intersection of each of these curves with the circle $\T$ has zero length. 
\par{\rm (b)} One has $\eta<|\th|<\eps$ on $\Ga_\eps\cap\D$ for some positive $\eta=\eta(\eps)$. 
\par{\rm (c)} The arc length $|dz|$ on $\Ga_\eps\cap\D$ is a Carleson measure, i.e., $H^1\subset L^1(\Ga_\eps,|dz|)$; moreover, the norm of the corresponding embedding operator is bounded by a constant depending only on $\eps$. 
\par{\rm (d)} For every $F\in H^1$, the equality 
$$\int_\T\f F\th dz =\int_{\Ga_\eps}\f F\th dz$$
holds, provided that the curves in the family $\Ga_\eps$ are oriented appropriately. 
\end{lem}

\par The following lemma reflects the well-known fact that, for $1<r<\infty$, one can identify $K^r_\th$ with $\left(K^{r'}_\th\right)^*$, the dual of the $K^{r'}_\th$ space with $r':=r/(r-1)$; a proof can be found in \cite{C1}. This is intimately related to the direct sum decomposition 
$$H^r=\th H^r\oplus K^r_\th,\qquad 1<r<\infty,$$
which in turn is a consequence of the M. Riesz theorem on conjugate functions (see \cite[Chapter III]{G}).

\begin{lem}\label{lem:cohnduality} Given $r\in(1,\infty)$, an inner function $\th$ and an element $f$ of $K^1_\th$, we have $f\in H^r$ (and hence $f\in K^r_\th$) if and only if $f$ induces a continuous linear functional on $K^{r'}_\th$, where $1/r+1/r'=1$, under the pairing 
$$\langle f,g\rangle=\int_\T\ov fg\,dm,\qquad g\in K^{r'}_\th.$$
Moreover, the norm $\|f\|_r$ is then comparable to the norm of the associated functional in $\left(K^{r'}_\th\right)^*$, with comparison constants depending only on $r$. 
\end{lem}

The next fact was established by Cohn in \cite{C1}. 

\begin{lem}\label{lem:cohnequivnorms} Given an inner function $\th$ and a number $\eps\in(0,1)$, let $\Ga_\eps=\Ga_\eps(\th)$ be the system of Carleson curves coming from Lemma \ref{lem:carlesoncurves}. If $f\in K^1_\th$ and $1<r<\infty$, then 
$$\|f\|_r\asymp\left(\int_{\Ga_\eps}|f(z)|^r|dz|\right)^{1/r},$$
with the understanding that $\|f\|_r=\infty$ whenever $f\notin L^r(\T,m)$. The constants involved in the $\asymp$ relation depend only on $r$ and $\eps$. 
\end{lem}
	
The following \lq\lq maximum principle", proved by Cohn in \cite{C2}, can be viewed as an endpoint version of the preceding result. 

\begin{lem}\label{lem:cohnmaxprinciple} Let $\th$ be an inner function, and let $f\in K^1_\th$. If $f$ is bounded on the set 
\begin{equation}\label{eqn:defsublevel}
\Om(\th,\eps):=\{z\in\D:\,|\th(z)|<\eps\}
\end{equation}
with some $\eps\in(0,1)$, then $f\in H^\infty$ and 
$$\|f\|_\infty\lesssim\sup\{|f(z)|:\,z\in\Om(\th,\eps)\},$$
where the comparison constant depends only on $\eps$. 
\end{lem}

The next result, borrowed from \cite{DSpb93} and/or \cite{DAJM}, characterizes the pairs $(f,\th)$, with $f\in\Al$ and $\th$ inner, where $f$ admits multiplication or division by every power of $\th$ in $\La^\al$. (Recall that this situation is also described by saying that $f$ is $(\al,\th)$-nice.) Once again, the sublevel set \eqref{eqn:defsublevel} will play a role in the statement. 

\begin{lem}\label{lem:kdmultdiv}
Suppose that $0<\al<\infty$, $n\in\N$, and $n>\al$. Given $f\in\Al$ and an 
inner function $\th$, the following conditions are equivalent. 
\par{\rm (i)} $f\th^n\in\Al$.
\par{\rm (ii)} $f\ov\th^n\in\La^{\al}$.
\par{\rm (iii)} $f\th^k\in\La^{\al}$ for all $k\in\Z$.
\par{\rm (iv)} For some (or each) $\eps\in(0,1)$, 
\begin{equation}\label{eqn:supsublevel}
\sup\left\{\f{|f(z)|}{(1-|z|)^\al}:\,z\in\Omte\right\}<\infty.
\end{equation}
\end{lem}

\par As usual, the equivalence relations above are actually accompanied by the appropriate estimates of the underlying quantities in terms of one another. In particular, letting $S=S(f,\th,\al,\eps)$ denote the supremum in \eqref{eqn:supsublevel} and recalling the notation 
$$\|f\|_{\al,\th}:=\|f\|_{\La^\al}+\|f\ov\th^n\|_{\La^\al},$$
with $n=[\al]+1$, we have 
\begin{equation}\label{eqn:sususu}
S\lesssim\|f\|_{\al,\th},
\end{equation}
where the constant involved is independent of $f$ and $\th$. 
\par We mention in passing that conditions (i) and (ii) in Lemma \ref{lem:kdmultdiv} cannot, in general, be replaced by $f\th\in\Al$ and $f\ov\th\in\La^{\al}$ when $\al\ge1$; see \cite[Chapter I]{Shi} and \cite{DScand, DSpb10} for more details. 
\par Finally, we supplement Lemma \ref{lem:kdmultdiv} with a related result from \cite[Section 4]{DSpb93}. 

\begin{lem}\label{lem:projlip} Suppose that $f\in H^1$ and $B$ is an interpolating Blaschke product with zeros $\{z_j\}$. For $0<\al<\infty$, one has 
$P_-(f\ov B)\in\La^\al$ if and only if 
$$\sup_j\f{|f(z_j)|}{(1-|z_j|)^\al}<\infty.$$
\end{lem}

\par Strictly speaking, the assumption in \cite{DSpb93} was that $f\in H^2$ rather than $f\in H^1$, but the proof given there works in the latter case as well.

\section{Proof of Theorem \ref{thm:factalphapr}} 

We note that $r\in(1,\infty)$ and define the conjugate exponent $r'$ by $1/r+1/r'=1$. In view of the duality 
between $K^r_\th$ and $K^{r'}_\th$ (recall Lemma \ref{lem:cohnduality}), the membership of $f$ in $K^r_\th$ (and hence in $H^r$) 
will be established as soon as we check that $f$ generates a continuous linear functional on $K^{r'}_\th$. The 
bound on the functional's norm, to be obtained along the way, will then yield the required estimate on $\|f\|_r$. 
The mapping $\ph\mapsto\int\ov f\ph$ is well defined at least for $\ph\in K^\infty_\th$; and since $K^\infty_\th$ 
is dense in $K^{r'}_\th$, it will suffice to show that the desired norm inequality holds on the smaller subspace. 
\par Now, for $\ph\in K^\infty_\th$, we have 
\begin{equation}\label{eqn:linfunfphi}
\int_\T f\ov\ph\,dm=\int_\T gh\ov\ph\,dm=\int_\T g\ov\ph\th\cdot h\ov\th\,dm=\int_\T\psi h\ov\th\,dm, 
\end{equation}
where we put $\psi:=g\ov\ph\th$. Since $g\in H^p$ and $\ov\ph\th\in H^\infty_0$, it follows that $\psi\in H^p_0$ 
and 
\begin{equation}\label{eqn:insertpminus}
\int_\T\psi h\ov\th\,dm=\int_\T\psi P_-\left(h\ov\th\right)\,dm, 
\end{equation}
this last integral being understood as 
\begin{equation}\label{eqn:radlim}
\lim_{\rho\to1^-}\int_\T\psi(\rho\ze)\cdot P_-\left(h\ov\th\right)(\ze)\,dm(\ze).
\end{equation}
\par To keep on the safe side, we remark that the function $\psi$, and hence the integrand 
on the right-hand side of \eqref{eqn:insertpminus}, need not be integrable in the case $0<p<1$, 
while the integral in \eqref{eqn:radlim} certainly makes sense. The existence of 
the limit in \eqref{eqn:radlim} is then due to Lemma \ref{lem:drsduality}, since $\al=p^{-1}-r^{-1}>p^{-1}-1$. 
Also, if one replaces $P_-$ by $P_+$ in \eqref{eqn:radlim}, the resulting quantity will be $0$; 
this explains the identity \eqref{eqn:insertpminus}. 
\par Next, we set $q:=(\al+1)^{-1}$ (so that $0<q<1$) and observe that 
\begin{equation}\label{eqn:pqrprime}
\f1q=\al+1=\f1p-\f1r+1=\f1p+\f1{r'}.
\end{equation}
Combining \eqref{eqn:linfunfphi} with \eqref{eqn:insertpminus} and using Lemma \ref{lem:drsduality} again, we find that 
\begin{equation}\label{eqn:shishka}
\left|\int_\T f\ov\ph\,dm\right|=\left|\int_\T\psi P_-\left(h\ov\th\right)\,dm\right|
\le c_\al\|\psi\|_q\left\|P_-\left(h\ov\th\right)\right\|_{\La^\al}, 
\end{equation}
with a suitable constant $c_\al>0$. Recalling \eqref{eqn:pqrprime} and the fact that $|\psi|=|g|\cdot|\ph|$, 
we now apply H\"older's inequality to get 
$$\|\psi\|_q\le\|g\|_p\|\ph\|_{r'}$$
and plug this into the rightmost side of \eqref{eqn:shishka}. This yields 
$$\left|\int_\T f\ov\ph\,dm\right|\le c_\al\|g\|_p\left\|P_-\left(h\ov\th\right)\right\|_{\La^\al}\|\ph\|_{r'}.$$ 
The function $\ph\in K^\infty_\th$ having been arbitrary, we conclude that the linear functional 
$\ph\mapsto\int\ov f\ph$ is indeed continuous on $K^{r'}_\th$, with norm at most 
$$\const\cdot\|g\|_p\left\|P_-\left(h\ov\th\right)\right\|_{\La^\al}.$$ 
Equivalently, $f$ is in $K^r_\th$ (and hence in $H^r$) and satisfies \eqref{eqn:normalphapr}, as required.

\section{Proofs of Corollaries \ref{cor:innsomealpha} and \ref{cor:smallinnfact}}

\medskip\noindent{\it Proof of Corollary \ref{cor:innsomealpha}.} Let 
$$E_f:=\{p\in(0,\infty):\,f\in H^p\}$$
and $p_0:=\sup E_f$. We want to prove that $E_f=(0,\infty)$, or equivalently, that  $p_0=\infty$. Assuming the contrary, while recalling that $f\in H^1$, we would have $1\le p_0<\infty$. We could then find a number $\al\in(0,p_0^{-1})$ such that $P_-(I\ov\th)\in\La^\al$. (Indeed, by hypothesis, the latter holds for all sufficiently small $\al>0$.) This done, we could furthermore take $p\in E_f$ close enough to $p_0$ so as to satisfy 
$$\f p{1-\al p}\left(=:r\right)>p_0.$$
Now, since $f\in H^p$, we also have $F\in H^p$, and Theorem \ref{thm:factalphapr} would tell us (when applied to $g=F$ and $h=I$) that $f\in H^r$. This, however, is incompatible with the fact that $r>p_0$, as the definition of $p_0$ shows. The contradiction convinces us that we actually have $p_0=\infty$.\quad\qed

\medskip\noindent{\it Proof of Corollary \ref{cor:smallinnfact}.} As we know from Lemma \ref{lem:projlip}, condition \eqref{eqn:smallatzeros} is equivalent to saying that $P_-(I\ov B)\in\La^\al$. Consequently, Corollary \ref{cor:innsomealpha} applies (with $\th=B$) and tells us that $f$ is in every $H^p$ with $0<p<\infty$. The same is then true for $F$, and so 
\begin{equation}\label{eqn:ktoto}
|F(z_j)|\lesssim(1-|z_j|)^{-\de},\qquad j=1,2,\dots,
\end{equation}
for every $\de>0$. We may now combine \eqref{eqn:smallatzeros} and \eqref{eqn:ktoto} to infer that 
\begin{equation}\label{eqn:chtoto}
|f(z_j)|=|F(z_j)|\cdot|I(z_j)|\lesssim(1-|z_j|)^{\al-\de},\qquad j=1,2,\dots,
\end{equation}
for any fixed $\de\in(0,\al)$. This last estimate implies, by virtue of Lemma \ref{lem:projlip}, that 
\begin{equation}\label{eqn:gdeto}
f\ov B=P_-(f\ov B)\in\La^{\al-\de}
\end{equation}
(we have also used the hypothesis that $f\in K^1_B$). 
\par In particular, we see from \eqref{eqn:gdeto} that $f\ov B\in L^\infty(\T)$, whence $f\in H^\infty$ and $F\in H^\infty$. The latter means that \eqref{eqn:ktoto}
actually holds with $\de=0$ as well. Proceeding as above, we then arrive at the improved versions of \eqref{eqn:chtoto} and \eqref{eqn:gdeto}, where the improvement consists in putting $\de=0$ throughout. Eventually we conclude that $f\ov B\in\La^\al$, or equivalently, that $\widetilde f\in\Al$. Because dividing a function by its inner factor preserves membership in $\Al$ (see, e.g., \cite{H} for a proof), the outer factor of $\widetilde f$ (which is $F$) will also be in $\Al$. \quad\qed

\section{Proof of Theorem \ref{thm:factalphaprbis}}

First let us assume that $\al<1/p<\al+1$, so that $1<r<\infty$. To deal with this case, we fix a number $\eps\in(0,1)$ and consider the Carleson curves $\Ga_\eps=\Ga_\eps(\th)$, as described in Lemma \ref{lem:carlesoncurves}. In view of Lemma \ref{lem:cohnequivnorms}, it is enough to show that the quantity $\int_{\Ga_\eps}|f|^r|dz|$ is finite and admits the appropriate bound. We have 
\begin{equation}\label{eqn:lrnormcurve}
\begin{aligned}
\int_{\Ga_\eps}|f(z)|^r|dz|&=\int_{\Ga_\eps}|g(z)|^r|h(z)|^r|dz|\\
&=\int_{\Ga_\eps}|g(z)|^p|g(z)|^{r-p}|h(z)|^r|dz|.
\end{aligned}
\end{equation}
Using the inequalities 
$$|g(z)|\le\|g\|_{A^{-1/p}}\cdot(1-|z|)^{-1/p},\qquad z\in\D,$$
and
\begin{equation}\label{eqn:huhuhu}
|h(z)|\lesssim\|h\|_{\al,\th}\cdot(1-|z|)^\al,\qquad z\in\Ga_\eps\cap\D,
\end{equation}
in conjunction with the identity $p^{-1}-r^{-1}=\al$, we find that 
$$|g(z)|^{r-p}|h(z)|^r\lesssim\|g\|_{A^{-1/p}}^{r-p}\|h\|_{\al,\th}^r,
\qquad z\in\Ga_\eps\cap\D.$$ 
(To verify \eqref{eqn:huhuhu}, recall Lemma \ref{lem:kdmultdiv} and the estimate \eqref{eqn:sususu} stated next to it.) It now follows from \eqref{eqn:lrnormcurve} that 
\begin{equation}\label{eqn:contcurve}
\int_{\Ga_\eps}|f(z)|^r|dz|\lesssim\|g\|_{A^{-1/p}}^{r-p}\|h\|_{\al,\th}^r
\int_{\Ga_\eps}|g(z)|^p|dz|.
\end{equation}
Recalling that arc length $|dz|\big|_{\Ga_\eps}$ is a Carleson measure, we get 
$$\int_{\Ga_\eps}|g(z)|^p|dz|\lesssim\|g\|_p^p,$$
where the constant involved depends only on $\eps$. Combining this last inequality with \eqref{eqn:contcurve} yields 
$$\int_{\Ga_\eps}|f(z)|^r|dz|\lesssim\|g\|_{A^{-1/p}}^{r-p}
\|g\|_p^p\|h\|_{\al,\th}^r,$$
or equivalently, 
\begin{equation}\label{eqn:endcurve}
\left(\int_{\Ga_\eps}|f(z)|^r|dz|\right)^{1/r}\lesssim\|g\|_{A^{-1/p}}^{\alpha p}
\|g\|_p^{1-\alpha p}\|h\|_{\al,\th}.
\end{equation}
Finally, since 
$$\left(\int_{\Ga_\eps}|f(z)|^r|dz|\right)^{1/r}\gtrsim\|f\|_r$$
by Lemma \ref{lem:cohnequivnorms}, we deduce from \eqref{eqn:endcurve} that 
$$\|f\|_r\lesssim\|g\|_{A^{-1/p}}^{\alpha p}\|g\|_p^{1-\alpha p}\|h\|_{\al,\th},$$
as desired. Moreover, the constant hidden in the $\lesssim$ sign can be taken to depend on $p$ and $\al$ only. (So far, it was also dependent on $\eps$, but we could have fixed the value of $\eps$, say $\eps=\f12$, from the start.) This proves 
\eqref{eqn:normalphaprbis} under the current assumption on the parameters. 
\par Now suppose that $\al=1/p$, so that $r=\infty$. To handle this case, we fix (once again) a number $\eps\in(0,1)$ and consider the corresponding sublevel set 
$$\Om(\th,\eps):=\{z\in\D:\,|\th(z)|<\eps\}.$$
Since 
$$|g(z)|\le\|g\|_{A^{-\al}}\cdot(1-|z|)^{-\al},\qquad z\in\D,$$
and
$$|h(z)|\lesssim\|h\|_{\al,\th}\cdot(1-|z|)^\al,\qquad z\in\Om(\th,\eps),$$
we infer that 
$$|f(z)|=|g(z)|\cdot|h(z)|\lesssim\|g\|_{A^{-\al}}\|h\|_{\al,\th},
\qquad z\in\Om(\th,\eps).$$
We may then invoke Lemma \ref{lem:cohnmaxprinciple} to conclude that $f\in H^\infty$ and \eqref{eqn:normalphainfty} holds true.

\section{Proof of Proposition \ref{prop:sharpness}} 

Let $\th$ be the Blaschke product with zeros $z_n:=1-2^{-n}$ ($n=1,2,\dots$). We note, first, that $\{z_n\}$ is an interpolating sequence for $H^\infty$ (cf. \cite[Chapter VII]{G}), and secondly, that $\th'\in H^\ga$ for each $\ga\in(0,1)$. The latter is due to a theorem of Protas \cite{P}; see also \cite{AC, GN}. In fact, a result from \cite{CN} tells us that $\th'$ actually belongs to the weak $H^1$ space, which implies the preceding property as well. 

\par Next, assuming that the numbers $\al$, $p$ and $\si$ satisfy the hypotheses of the current proposition, we define 
the \lq\lq critical exponent" $r$ as in Theorem \ref{thm:factalphapr}, so that 
$$r:=\f p{1-\al p}>1,$$ 
and then fix a number $s$ with $r<s<\si$. This done, we can find a function $f\in K^s_\th\setminus H^\si$, a fact 
which we now explain. Indeed, let $\{w_n\}$ be a sequence of complex numbers such that 
\begin{equation}\label{eqn:sconverges}
\sum_n|w_n|^s(1-|z_n|)<\infty,
\end{equation}
but 
\begin{equation}\label{eqn:sigmadiverges}
\sum_n|w_n|^\si(1-|z_n|)=\infty
\end{equation}
(one possible choice is $w_n=n^{-1}(1-|z_n|)^{-1/s}$). Further, let $\Phi\in H^s$ be a solution to the interpolation 
problem
\begin{equation}\label{eqn:interprob}
\Phi(z_n)=w_n,\qquad n=1,2,\dots,
\end{equation}
the existence of such a $\Phi$ being guaranteed by \eqref{eqn:sconverges} in view of the Shapiro--Shields 
theorem (see \cite{SS} or \cite[Chapter VII]{G}). Finally, we put 
\begin{equation}\label{eqn:decomp}
f:=\Phi-\th P_+(\ov\th\Phi)
\end{equation}
and note that $f\in K^s_\th(\subset K^1_\th)$; in fact, $f$ is the projection of $\Phi$ onto $K^s_\th$ parallel 
to $\th H^s$. (The direct sum decomposition 
$$H^s=K^s_\th\oplus\th H^s$$ 
holds by the M. Riesz theorem, since $1<s<\infty$.) It follows from \eqref{eqn:interprob} and 
\eqref{eqn:decomp} that $f(z_n)=w_n$ for all $n$. Consequently, the quantity $\sum_n|w_n|^\si(1-|z_n|)$ can be 
written as $\int|f|^\si d\mu$, where $\mu$ is the discrete (Carleson) measure with weights $1-|z_n|$ at 
the $z_n$'s; and if $f$ were in $H^\si$, this quantity would have to be finite, in contradiction 
with \eqref{eqn:sigmadiverges}. We conclude that $f\in K^s_\th\setminus H^\si$, as desired. 
\par Now, we obviously have $f=gh$, where 
$$g(z):=f(z)\cdot(1-z)^{-\al}$$
and 
$$h(z):=(1-z)^\al.$$
The proof will be complete as soon as we check that $g\in H^p$ and $h$ is $(\al,\th)$-nice. To this end, we first 
recall that $s>r>p$ and let $q$ be the exponent defined by 
\begin{equation}\label{eqn:pqs}
q^{-1}+s^{-1}=p^{-1}.
\end{equation}
Since $s^{-1}<r^{-1}=p^{-1}-\al$, it follows that $q^{-1}>\al$ and the function $1/h(z)=(1-z)^{-\al}$ is, 
therefore, in $H^q$. Combining this with the fact that $f\in H^s$ and bearing \eqref{eqn:pqs} in mind, 
we see that $g(=f/h)$ is in $H^p$, thanks to H\"older's inequality.
\par Finally, we take $\eps\in(0,1)$ to be appropriately small and observe that the sublevel set $\Om(\th,\eps)$, defined as in \eqref{eqn:defsublevel}, lies 
in a certain Stolz angle with vertex at $1$, i.e., in the convex hull of $t\D\cup\{1\}$ with some $t\in(0,1)$. (Indeed, since $\th$ is an interpolating Blaschke product, $\Om(\th,\eps)$ is contained in the union of hyperbolic disks of some fixed radius with centers at the $z_n$'s; 
see \cite[Chapter X, Lemma 1.4]{G}.) Because $h\in\Al$ and 
$$|h(z)|=|1-z|^\al\lesssim(1-|z|)^\al$$
on such a Stolz angle, the desired conclusion that $h$ is $(\al,\th)$-nice is a consequence of Lemma \ref{lem:kdmultdiv}.

\section{Proof of Theorem \ref{thm:embedlip}}

{\rm (A)} Assuming that $\al>p^{-1}$, put $\be=\al-p^{-1}$ and $s=(\be+1)^{-1}$. 
To prove that $\widetilde f:=\ov z\ov f\th$ lies in $A^\be$ and satisfies the norm estimate \eqref{eqn:normalphabeta}, we invoke Lemma \ref{lem:drsduality} which allows us to identify $A^\be$ with the dual of $H^s$. Thus, we are going to show that $\widetilde f$ generates a continuous linear functional on $H^s$ whose norm admits the appropriate bound. Moreover, since $H^\infty$ is a dense subspace of $H^s$, it will suffice to look at the action of the corresponding functional 
\begin{equation}\label{eqn:phifunct}
\ph\mapsto\int_\T\ov{\widetilde f}\ph\,dm
\end{equation}
on a function $\ph\in H^\infty$. For such a $\ph$, we have
\begin{equation}\label{eqn:palka}
\int_\T\ov{\widetilde f}\ph\,dm=\int_\T zf\ov{\th}\ph\,dm=
\int_\T z\ph g\cdot h\ov{\th}\,dm.
\end{equation}
The last quantity equals 
$$\int_\T z\ph g\cdot P_-(h\ov{\th})\,dm:=\lim_{\rho\to1^-}
\int_\T\ze\ph(\rho\ze)g(\rho\ze)\cdot P_-(h\ov{\th})(\ze)\,dm(\ze)$$
(where the limit exists by virtue of Lemma \ref{lem:drsduality}), because the contribution of $P_+(h\ov{\th})$ to the integral is zero. We may therefore rewrite \eqref{eqn:palka} as
\begin{equation}\label{eqn:kopalka}
\int_\T\ov{\widetilde f}\ph\,dm=\int_\T z\ph g\cdot P_-(h\ov{\th})\,dm.
\end{equation}
Using Lemma \ref{lem:drsduality} once again, we infer from \eqref{eqn:kopalka} that 
\begin{equation}\label{eqn:skalka}
\left|\int_\T\ov{\widetilde f}\ph\,dm\right|\lesssim
\|\ph g\|_q\left\|P_-(h\ov{\th})\right\|_{\La^\al},
\end{equation}
where $q=(\al+1)^{-1}$. Furthermore, since 
\begin{equation}\label{eqn:tselka}
\f1q=\al+1=\f1p+\be+1=\f1p+\f1s,
\end{equation}
H\"older's inequality gives 
$$\|\ph g\|_q\le\|g\|_p\|\ph\|_s,$$
and we plug this into the right-hand side of \eqref{eqn:skalka} to get 
\begin{equation}\label{eqn:microzhopa}
\left|\int_\T\ov{\widetilde f}\ph\,dm\right|\lesssim
\|g\|_p\left\|P_-(h\ov{\th})\right\|_{\La^\al}\|\ph\|_s.
\end{equation}
\par Hence we conclude that the functional \eqref{eqn:phifunct}, defined initially for $\ph\in H^\infty$, extends continuously to $H^s$ with norm at most 
\begin{equation}\label{eqn:samovar}
\const\cdot\|g\|_p\left\|P_-(h\ov{\th})\right\|_{\La^\al}.
\end{equation}
In view of Lemma \ref{lem:drsduality}, this proves the inclusion $\widetilde f\in A^\be$ and the norm estimate \eqref{eqn:normalphabeta}. 

\smallskip{\rm (B)} Assuming that $\al=p^{-1}$, we define the exponents $\be$ and $s$ as before, which yields $\be=0$ and $s=1$. The rest of the proof is almost identical to what we did in part (A), provided that we agree to interpret $A^0$ (i.e., the $A^\be$ space with $\be=0$) as $\bmoa$. More precisely, proceeding as above, we arrive at \eqref{eqn:microzhopa}, which now reads 
$$\left|\int_\T\ov{\widetilde f}\ph\,dm\right|\lesssim
\|g\|_p\left\|P_-(h\ov{\th})\right\|_{\La^\al}\|\ph\|_1.$$
This means that the functional \eqref{eqn:phifunct}, defined initially on $H^\infty$, extends continuously to $H^1$ with norm bounded by \eqref{eqn:samovar}. The ($H^1$, $\bmoa$) duality theorem (see, e.g., \cite[Chapter VI]{G}) finally allows us to conclude that $\widetilde f$ is in $\bmoa$ and satisfies \eqref{eqn:normalphabmo}.

\section{Proof of Proposition \ref{prop:sharpnesslip}}

As in the proof of Proposition \ref{prop:sharpness}, we take $\th$ to be the Blaschke product with zeros $z_j=1-2^{-j}$ ($j=1,2,\dots$). We know that $\{z_j\}$ is an interpolating sequence and that $\th'\in\bigcap_{0<\ga<1}H^\ga$. 
\par Next, we put $\be=\al-p^{-1}$ and pick a number $\eps>0$ such that 
$$0<\eps<\min\left(\de-\be,\,p^{-1}\right).$$ 
Now let $\Phi\in H^\infty$ be a function that solves the interpolation problem 
$$\Phi(z_j)=(1-|z_j|)^{\be+\eps}\qquad(j=1,2,\dots),$$
and let $f$ be the orthogonal projection (in $H^2$) of $\Phi$ onto $K^2_\th$. Thus, 
$$f=\Phi-\th P_+(\ov\th\Phi)=\th P_-(\ov\th\Phi).$$
Then $f(z_j)=\Phi(z_j)$ for each $j$, so that 
\begin{equation}\label{eqn:valuesoff}
f(z_j)=(1-|z_j|)^{\be+\eps}\qquad(j=1,2,\dots),
\end{equation}
and Lemma \ref{lem:projlip} tells us that $P_-(f\ov\th)\in\La^{\be+\eps}$. Because $f\in K^2_\th$, we have 
\begin{equation}\label{eqn:zftilde}
\ov{\widetilde f}=zf\ov\th=zP_-(f\ov\th),
\end{equation}
and it follows that $\widetilde f\in A^{\be+\eps}$. In particular, $\widetilde f$ is bounded, and so is $f$. 
\par Furthermore, another application of Lemma \ref{lem:projlip} in conjunction with \eqref{eqn:valuesoff} and \eqref{eqn:zftilde} shows that $\widetilde f\notin A^\de$. Indeed, since $\be+\eps<\de$, the values $f(z_j)$ tend to $0$ at an essentially slower rate than $(1-|z_j|)^\de$, so the lemma ensures that the last (and hence also the first) member of identity \eqref{eqn:zftilde} is not in $\La^\de$. 
\par Finally, we let $\tau=p^{-1}-\eps$ and write $f=gh$, where 
$$g(z):=(1-z)^{-\tau}\quad\text{\rm and}\quad h(z):=(1-z)^\tau f(z).$$
We have then $g\in H^p$ (because $\tau p<1$) and $h\in H^\infty$ (because $f\in H^\infty$). It is also true that $P_-(h\ov\th)\in\La^\al$. To verify this last claim, observe that 
$$|h(z_j)|=|1-z_j|^\tau|f(z_j)|=(1-|z_j|)^{\tau+\be+\eps}=(1-|z_j|)^\al$$
and invoke Lemma \ref{lem:projlip} once again. This done, we conclude that the current factors $g$ and $h$ have all the required properties. The proof is therefore complete.

\section{Proof of Theorem \ref{thm:embedlipbis}}

As in the proof of Theorem \ref{thm:embedlip}, we put $\be=\al-p^{-1}$ and $s=(\be+1)^{-1}$. Once again, our plan is to exploit the duality between $H^s$ and $A^\be$. More precisely, we shall arrive at \eqref{eqn:normalphabetabis} by showing that the functional \eqref{eqn:phifunct}, defined initially on $H^\infty$, extends continuously to $H^s$ and admits the appropriate norm estimate. 
\par For $\ph\in H^\infty$, we have \eqref{eqn:palka} as before, and we rewrite the resulting identity as 
\begin{equation}\label{eqn:violin}
\int_\T\ov{\widetilde f}\ph\,dm=
\f1{2\pi i}\int_\T\f{\ph gh}{\th}\,dz.
\end{equation}
Furthermore, the last integral remains unchanged if the integration set, $\T$, gets replaced by the (appropriately oriented) Carleson curves $\Ga_\eps=\Ga_\eps(\th)$ with some, or any, $\eps\in(0,1)$. This is guaranteed by property (d) from Lemma \ref{lem:carlesoncurves}. Consequently, once $\eps$ is fixed, \eqref{eqn:violin} implies that 
\begin{equation}\label{eqn:cello}
\begin{aligned}
\left|\int_\T\ov{\widetilde f}\ph\,dm\right|
&=\left|\f1{2\pi i}\int_{\Ga_\eps}\f{\ph gh}{\th}\,dz\right|\\
&\le\f1{2\pi}\int_{\Ga_\eps}\f{|\ph gh|}{|\th|}\,|dz|
=\f1{2\pi}\int_{\Ga_\eps}\f{|u|}{|\th|}\,|dz|,
\end{aligned}
\end{equation}
where we write $u:=\ph gh$. The integrals over $\Ga_\eps$ are actually taken over $\Ga_\eps\cap\D$, since $m(\Ga_\eps\cap\T)=0$. 
\par Now, for $z\in\D$, we have 
$$|\ph(z)|\le\|\ph\|_s\cdot(1-|z|)^{-1/s},$$
whence 
\begin{equation*}
\begin{aligned}
|\ph(z)|&=|\ph(z)|^{1-s}|\ph(z)|^s\\
&\le\|\ph\|^{1-s}_s\cdot(1-|z|)^{-(1-s)/s}\cdot|\ph(z)|^s,
\end{aligned}
\end{equation*}
or equivalently, 
\begin{equation}\label{eqn:viola}
|\ph(z)|\le\|\ph\|^{1-s}_s\cdot(1-|z|)^{-\be}\cdot|\ph(z)|^s.
\end{equation}
Also, 
\begin{equation}\label{eqn:bassoon}
|g(z)|\le\|g\|_{A^{-1/p}}\cdot(1-|z|)^{-1/p},\qquad z\in\D,
\end{equation}
and 
\begin{equation}\label{eqn:piano}
|h(z)|\lesssim\|h\|_{\al,\th}\cdot(1-|z|)^\al,\qquad z\in\Ga_\eps\cap\D.
\end{equation}
(To verify the latter estimate, recall Lemma \ref{lem:kdmultdiv} and the discussion following it.) 
\par Multiplying \eqref{eqn:viola}, \eqref{eqn:bassoon} and \eqref{eqn:piano} together, while keeping the relation $\be+p^{-1}=\al$ in mind, we find that 
$$|u(z)|\lesssim\|g\|_{A^{-1/p}}\|h\|_{\al,\th}\|\ph\|^{1-s}_s\cdot|\ph(z)|^s,
\qquad z\in\Ga_\eps\cap\D.$$ 
We now couple this with the inequality 
$$|\th(z)|\ge\eta(\eps),\qquad z\in\Ga_\eps\cap\D,$$
to estimate the last (rightmost) integral in \eqref{eqn:cello}. This yields 
\begin{equation}\label{eqn:medved}
\begin{aligned}
\int_{\Ga_\eps}\f{|u|}{|\th|}\,|dz|&\lesssim
\|g\|_{A^{-1/p}}\|h\|_{\al,\th}\|\ph\|^{1-s}_s\int_{\Ga_\eps}|\ph(z)|^s\,|dz|\\
&\lesssim\|g\|_{A^{-1/p}}\|h\|_{\al,\th}\|\ph\|_s,
\end{aligned}
\end{equation}
where the last step relies on the fact that the arc length $|dz|\big|_{\Ga_\eps\cap\D}$ is a Carleson measure. 
\par Finally, we combine \eqref{eqn:cello} and \eqref{eqn:medved} to obtain  
$$\left|\int_\T\ov{\widetilde f}\ph\,dm\right|\lesssim
\|g\|_{A^{-1/p}}\|h\|_{\al,\th}\|\ph\|_s.$$ 
This means that the functional \eqref{eqn:phifunct} admits a continuous extension to $H^s$, with norm not exceeding 
$$\const\cdot\|g\|_{A^{-1/p}}\|h\|_{\al,\th}.$$
Moreover, by assigning a concrete numerical value to $\eps$, we may assume that the constant is independent of that parameter (and therefore depends on $\al$ and $p$ only). An application of Lemma \ref{lem:drsduality} now leads us to the desired conclusion that $\widetilde f$ lies in $A^\be$ and satisfies \eqref{eqn:normalphabetabis}.

\end{document}